\documentclass [oneside,a4paper,11pt] {article}
\usepackage{latexsym,epsfig,graphpap,graphics,amsmath,amssymb,}
\begin{document}

\begin{center}
{\bf{PARALLEL VECTOR FIELDS ON THE NONINVARIANT HYPERSURFACE OF A SASAKIAN MANIFOLD}}\\
{\it By}

\vspace{.3cm} \noindent{\bf{Sachin Kumar Srivastava, Alok Kumar Srivastava
and Dhruwa Narain}}\\

\end{center}

{\bf Abstract~:} In 1970, Samuel I. Goldberg and Kentaro Yano defined the notion of noninvariant hypersurface of a Sasakian manifold [1]. In this paper we have studied the properties of parallel vector fields with respect to induced connection on the noninvariant hypersurface $M$ of a Sasakian manifold $\tilde M$ with $(\phi, g, u, v, \lambda)-$ structure and proved that if the vector field $V$ is parallel with respect to induced connection on $M$ then $M$ is totally geodesic.

\vspace{.3cm}\noindent{\bf{2000 Mathematics Subject Classification :}}~ 53D10, 53C25, 53C21.

\vspace{.5cm}
\noindent{\bf{1.~~Introduction~:}}~Let $\tilde M$ be a $(2n + 1)-$dimensional Sasakian manifold with a tensor field $\phi$ of type (1, 1),  a fundamental vector field $\xi$ and $1-form \, \eta$ such that\\
\nolinebreak
(1.1)\hspace{1in}$\eta(\xi) = 1$\\
(1.2)\hspace{1in}$\tilde \phi^2 = - I + \eta \otimes \xi $\\
\vspace{.2cm}
where $I$ denotes the identity transformation.\\
(1.3)(a)\hspace{.8in}
$\eta \, o\,\phi = 0\qquad (b) \qquad \phi \xi = 0 \qquad (c) \qquad rank(\tilde \phi) = 2n$\\
\vspace{.2cm}
If $\tilde M$ admits a Riemannian metric $\tilde g$, such that\\
\vspace{.2cm}
(1.4)\hspace{1in}$\tilde g(\tilde\phi X,\tilde\phi Y) = \tilde g(X, Y) - \eta(X) \eta(Y)$\\
\vspace{.2cm}
(1.5)\hspace{1in}$ \tilde g(X,\xi)=\eta(X)$\\
\vspace{.2cm}
then $\tilde M$ is said to admit a $(\tilde\phi,\xi,\eta,\tilde g)-$ structure called contact metric strucure.\\
If moreover,\\
(1.6)\hspace{1in}$(\tilde\nabla_X \tilde\phi )\, Y =\tilde g (X, \, Y) \xi - \eta (Y) X $\\
 \hspace{.3in} and\\
(1.7)\hspace{1in}$\tilde\nabla_X\xi=-\tilde\phi X$\\
where $\tilde\nabla$ denotes the Riemannian connection of the Riemannian metric g, then $(\tilde M,\tilde\phi,\xi,\eta,\tilde g)$ is called a Sasakian manifold [9].
\vspace{3mm}\parindent=8mm
If we define ${}^\prime F(X, Y) = g (\phi X, Y)$, then in addition to above relation we find\\
(1.8)\hspace{1in}${}^\prime F(X, Y) + {}^\prime F(Y, X) = 0 $\\
(1.9)\hspace{1in}${}^\prime F(X, \phi Y) = {}^\prime F(Y, \phi X)$\\
(1.10)\hspace{.92in}${}^\prime F(\phi X, \phi Y) = {}^\prime F(X, Y)$\\
\vspace{.5cm}
\noindent{\bf{2.~~Noninvariant hypersurface of a Sasakian manifold~:}}~\\
\vspace{.1cm}
\hspace{.5cm}Let us consider a $2n-$dimensional manifold $M$ embedded in $\tilde M$ with embedding $b : M \rightarrow \tilde M$. The map $b$ induces a linear transformation map $B$ (called Jacobian map), $B : T_p \rightarrow T_{b_p}$.
\vspace{.1cm}
Let an affine normal $N$ of $M$ is in such a way that $\tilde\phi N$ is always tangent to the hypersurface and satisfying the linear transformations

\vspace{.25cm}\noindent
(1.5) \hspace{1in} $\tilde\phi BX = B \phi X + u(X) N$\\
(1.6) \hspace{1in} $\tilde\phi N = - BU$\\
(1.7) \hspace{1in} $\xi = BV + \lambda N$\\
(1.8) \hspace{1in} $\eta(BX) = v(X)$

\vspace{.5cm}\noindent
where $\phi$ is a (1, 1) type tensor; $U,\, V$ are vector fields; $u, v$ are $1- form$ and $\lambda$ is a $C^\infty -$ function. If $u \ne 0$, $M$, is called a noninvariant hypersurface of $\tilde M$ [1].

\vspace{.3cm}\parindent=8mm
Operating (2.1), (2.2), (2.3) and (2.4)  by $\tilde\phi$ and using (1.1), (1.2) and (1.3) and taking tangent normal parts separately, we get the following induced structure on $M$,

\vspace{.3cm}\noindent
(2.5)(a) \hspace{.96in} $\phi^2 X = -  X + u (X) U + v (X) V$

\vspace{.3cm}\parindent=8mm
(b)\hspace{1in} $u(\phi X) = \lambda v (X), \qquad v(\phi X) =  -\, \eta (N) \, u(X)$

\vspace{.3cm}\parindent=8mm
(c)\hspace{1in} $\phi U = - \, \eta(N) \, V , \qquad \phi V = \lambda U$

\vspace{.3cm}\parindent=8mm
(d)\hspace{1in} $u(U) = 1 -  \lambda \eta(N), \qquad u(V) = 0$

\vspace{.3cm}\parindent=8mm
(e)\hspace{1in} $v(U) = 0, \qquad v(V) = 1 -   \lambda \eta(N)$

\vspace{.3cm}\noindent
and from (1.4) and (1.5), we get the induced metric $g$ on $M$, i.e.,

\vspace{.3cm}
\noindent
(2.6) \hspace{1in} $g(\phi X, \phi Y) = g (X, Y) - u(X) u(Y) - v(X) v(Y)$\\
(2.7) \hspace{1in} $g (U, X) = u (X), \qquad g (V, X) = v (X).$

\vspace{.3cm}
\parindent=8mm
If we consider $\eta(N) = \lambda$, we get the following structures on $M$.

\vspace{.3cm}
\noindent
(2.8)(a)\hspace{1.1in} $\phi^2 = -  I +  u \otimes U + v \otimes V$

\vspace{.3cm}
\parindent=8mm
(b) \hspace{1in} $\phi U = - \, \lambda V, \qquad \phi V = \lambda U$

\vspace{.3cm}
\parindent=8mm
(c) \hspace{1in} $u \, \circ \phi = \lambda v, \qquad v \, \circ \, \phi =  - \lambda u$

\vspace{.3cm}
\parindent=8mm
(d) \hspace{1in} $u( U) = 1 - \lambda^2 , \qquad u(V) = 0$

\vspace{.3cm}
\parindent=8mm
(e) \hspace{1in} $v (U) = 0,  \qquad v (V) = 1 - \lambda^2.$

\vspace{.3cm}
\parindent=8mm
A manifold $M$ with a metric $g$ satisfying (2.6), (2.7) and (2.8) is called manifold with $(\phi, g, u, v, \lambda )-$structure[2].
\vspace{.3cm}
\parindent=8mm
Let $\nabla$ be the induced connection on the hypersurface $M$ of the affine connection $\tilde\nabla$ of $\tilde M$.

\vspace{.3cm}
\parindent=8mm
Now using Gauss and Weingarten's equations

\vspace{.3cm}
\noindent
(2.9)\hspace{1in} $\tilde\nabla _{BX} BY = B\nabla_X Y + h (X, Y) N$

\vspace{.3cm}
\noindent
(2.10) \hspace{.9in} $\tilde\nabla_{BX} N = BHX + w (X) N,$ ~where~ $g (HY, Z) = h(Y, Z).$

\vspace{.3cm}
\parindent=8mm
Here $h$ and $H$ are the second fundamental tensors of type (0, 2) and (1, 1) and $w$ is a $1-form$. Now differentiating (2.1), (2.2), (2.3) and (2.4) covariantly and using (2.9), (2.10), (1.6) and reusing (2.1), (2.2), (2.3) and (2.4), we get

\vspace{.3cm}
\noindent
(2.11) \hspace{.9in} $(\nabla_Y \phi) (X) = v (X) Y - g (X, Y) V - h (X, Y) U - u (X) HY$

\vspace{.3cm}
\noindent
(2.12) \hspace{.9in} $(\nabla_Y u) (X) = - h(\phi X, Y) - u (X) w (Y) - \lambda g (X, Y)$

\vspace{.3cm}
\noindent
(2.13) \hspace{.9in} $(\nabla_Y v) (X) = g(\phi Y, X) + \lambda h (X, Y)$

\vspace{.3cm}
\noindent
(2.14) \hspace{.9in} $\nabla_Y U = w (Y) U - \phi HY - \lambda Y$

\vspace{.3cm}
\noindent
(2.15) \hspace{.9in} $\nabla_Y V = \phi Y + \lambda H Y$

\vspace{.3cm}
\noindent
(2.16) \hspace{.9in} $h(Y, V) = u(Y) - Y \lambda - \lambda w(Y)$

\vspace{.3cm}
\noindent
(2.17) \hspace{.9in} $h(Y, U) = - u(HY)$

\vspace{.3cm}
\parindent=8mm
Since $h (X, Y) = g (HX, Y)$, then from (1.5), and (2.17), we get

\vspace{.3cm}
\noindent
(2.18) \hspace{.9in} $h (Y, U) = 0 \Rightarrow HU = 0.$

\vspace{.5cm}
\noindent{\bf{3.~~Parallel vector fields with respect to induced connection~:}}~

\vspace{.3cm}
\parindent=8mm
Let $M$ be the noninvariant hypersurface with $(\phi, g, u, v, \lambda)-$structure of a Sasakian manifold $\tilde M$. 
{\it{A vector field $P$ is parallel with respect to the connection $\tilde\nabla$ if $\tilde\nabla_X P = 0, \forall X \in \Gamma (T\tilde M)$ }}[9].\\
\vspace{.25cm}\hspace{.2in}
If $\phi$ is a parallel vector field, then from (2.11) we have\\
(3.1) \hspace{.9in} $v(X)Y - g(X, Y)V + h(X, Y)U + u(X)HY = 0$\\
\noindent
Operating (3.1) by $u$, using (2.8)(a) and (2.18) we have\\
\vspace{.2cm}
\hspace{1.25in}$ (1 - \lambda^2)h(X, Y) = - u(Y)v(X)$\\
Putting $X = V$ in above equation we get $h(Y, V) = 0.$\\
Again operating (3.1) by $v$, using (2.8)(e) and $v(HY) = 0$, we get\\
\vspace{.2cm}\hspace{1.25in}$(1 - \lambda^2)g(X, Y) = v(X)v(Y)$\\
Using $h(Y, V) = 0$ in (2.16), we get $w = \frac{u}{\lambda} - d(log\lambda)$.\\
\vspace{.2cm}This leads to the following theorem :\\
{\bf {Theorem~3.1~:}}~{\it{Let $M$ be the noninvariant hypersurface of a Sasakian manifold $\tilde M$ with $(\phi, g, u, v, \lambda )-$ structure . If $\phi$ is parallel vector field on $M$ then, we have,}}$ for all X,y in \gamma (TM)$\\\\
\vspace{.2cm}(3.2)\hspace{1in}$(1 - \lambda^2)h(X, Y) = - u(Y)v(X)$\\
\vspace{.2cm}(3.3)\hspace{1in}$h(X, V) = 0$\\
\vspace{.2cm}(3.4)\hspace{1in}$(1 - \lambda^2)g(X, Y) = v(X)v(Y)$\\
\vspace{.2cm}(3.5)\hspace{1in}$w = \frac{u}{\lambda} - dlog\lambda$\\
\vspace{.2cm}
If $U$ is parallel vector field, then from (2.14) we have\\
\vspace{.2cm}\hspace{1.25in}$w(X)U - \phi HX - \lambda X = 0.$\\
\vspace{.2cm}Operating by $u$ in the both sides of the above equation and using (2.8)(d) and (2.16), we get\\
\vspace{.2cm}\hspace{1.25in}$h(X, \phi Y) = \lambda g(X, Y) - w(Y)u(X)$ and $w = 2\lambda u - \lambda^2 d(log\lambda)$.\\
\vspace{.2cm}This leads to the following theorem:\\
{\bf{Theorem~3.2~:}}~{\it{Let $M$ be the noninvariant hypersurface of a Sasakian manifold $\tilde M$ with $(\phi, g, u, v, \lambda )-$ structure . If $U$ is parallel vector field on $M$ then, we have}}\\
\vspace{.2cm}(3.6)\hspace{1in}$h(X, \phi Y) = \lambda g(X, Y) - w(Y)u(X)$\\
\vspace{.2cm}(3.7)\hspace{1in}$w = 2\lambda u - \lambda^2 d(log\lambda)$\\
\vspace{.2cm}If $V$ is parallel vector field, then from (2.16)we have\\
\vspace{.2cm}\hspace{1.25in}$ \phi Y + \lambda HY = 0 \implies \lambda h(X, Y) = g(\phi X, Y).$\\
Since $h$ is symmetric, above equation yields that $2\lambda h(X, Y) = 0$ where $\lambda \ne 0$, therefore $h = 0.$\\
\vspace{.2cm}This leads to the following theorem:\\
{\bf{Theorem~3.3~}}~{\it{Let $M$ be the noninvariant hypersurface of a Sasakian manifold $\tilde M$ with $(\phi, g, u, v, \lambda )-$ structure . If $V$ is parallel vector field on $M$ then, $M$ is totally geodesic.}}\\\\
\vspace{.2cm}If $h = 0$, then from (2.16) we have $w = \frac{u}{\lambda} - d(log\lambda) .$\\
\vspace{.2cm}This leads to the following theorem:\\
{\bf{Theorem~3.4~:}~}{\it {Let $M$ be the noninvariant hypersurface of a Sasakian manifold $\tilde M$ with $(\phi, g, u, v, \lambda )-$ structure . If either $V$ or $\phi$ or both are parallel vector field on $M$, then we must have}}\\
\vspace{.2cm}(3.8)\hspace{1.25in}$w = \frac{u}{\lambda} - d(log\lambda).$
\begin{center}
{\bf References}
\end{center}

\begin{enumerate}
\item[{[1]}] S. I. Goldberg and K. Yano : {\em Noninvariant hypersurfaces of almost contact manifolds}, J. Math. Soc., Japan, {\bf 22} (1970), 25-34. 
\item[{[2]}] D. Narain : {\em Hypersurfaces with $(f, g, u, v, \lambda)-$structure of an affinely cosymplectic manifold}, Indian Jour. Pure and Appl. Math.{\bf 20} No. 8 (1989), 799-803.
\item[{[3]}] Dhruwa Narain and S. K. Srivastva : {\em Noninvariant Hypersurfaces of Sasakian Space Forms}, Int. J. Contemp. Math. Sciences, Vol. {\bf 4}, 2009, no. 33, 1611 - 1617.
\item[{[4]}] Dhruwa Narain, S. K. Srivastava and Khushbu Srivastava : {\em A note of noninvariant  hypersurfaces of para Sasakian manifold}, IOSR Journal of Engineering(IOSRJEN) Vol.{\bf 2}, issue 2 (2012), 363-368.
\item[{[5]}] Dhruwa Narain and S. K. Srivastava : {\em On the hypersurfaces of almost r-contact manifold}, J. T. S. Vol.{\bf 2}(2008), 67-74.
\item[{[6]}] Rajendra Prasad and M. M. Tripathi : {\em Transversal Hypersurfaces of Kenmotsu Manifold}, Indian Jour. Pure and Appl. Math.{\bf 34} No. 3 (2003), 443-452.
\item[{[7]}] B. B. Sinha and Dhruwa Narain : {\em Hypersurfaces of nearly Sasakian manifold}, Ann. Fac. Sci. Dekishasa Zaire, {\bf 3} (2), (1977), 267-280.
\item[{[8]}] M. Okumura : {\em Totally umbilical hypersurface of a locally product Riemannian manifold}, Kodai Math. Sem. Rep.{\bf 19} (1967), 35-42.
\item[{[9]}] D. E. Blair  : {\em Almost contact manifolds with killing structure tensor}, Pacific J. of Math. {\bf 39} (12), (1971), 285-292.
\item[{[10]}] T. Miyazawa and S. Yamaguchi : {\em Some theorems on K-contact metric manifolds and Sasakian manifolds}, T. R. V. Math. {\bf 2} (1966), 40-52.
\end{enumerate}
\vspace{.2cm}
\hspace{.5cm}{\it{ Authors' addresses}:}\\
{Sachin Kumar Srivastava}\\
{Amity Institute of Applied Sciences,Amity University , Noida, U.P., India}\\
\vspace{.2cm}
$E-mail: sachink.ddumath@gmail.com$\\
Alok Kumar Srivastava\\
Department of Mathematics ,Govt. Degree College, Chunar -
Mirzapur,U.P., India \\
\vspace{.2cm}
$E-mail: aalok\_sri@yahoo.co.in$\\
Dhruwa Narain\\
Department of Mathematics and Statistics\\
D.D.U. Gorakhpur University Gorakhpur, Gorakhpur-273009, India\\
$E- mail:dhruwanarain\_dubey@yahoo.co.in$\\

\end{document}